\documentclass[11pt]{article}

\usepackage{IEEEtrantools}
\usepackage{amsfonts}
\usepackage{amssymb}
\usepackage{amsmath,amsfonts,amsthm}
\usepackage{mathrsfs}
\usepackage{url}
\usepackage{tikz}
\usetikzlibrary{matrix,arrows,decorations.pathmorphing}

\usepackage{color}
\newcommand{\bbbt}{\mathbb{T}}
\newcommand{\scrt}{\mathscr{T}}

\newcommand{\be}{\begin{equation}}
\newcommand{\ee}{\end{equation}}
\newcommand{\bea}{\begin{eqnarray}}
\newcommand{\eea}{\end{eqnarray}}
\newcommand{\bean}{\begin{eqnarray*}}
\newcommand{\eean}{\end{eqnarray*}}
\newcommand{\brray}{\begin{array}}
\newcommand{\erray}{\end{array}}
\newcommand{\biearray}{\begin{IEEEarray}{rCl}}
\newcommand{\eiearray}{\end{IEEEarray}}

\newcommand{\newsection}[1]{\setcounter{equation}{0}
\setcounter{dfn}{0}
\section{#1}}
\newtheorem{dfn}{Definition}[section]
\newtheorem{thm}[dfn]{Theorem}
\newtheorem{lmma}[dfn]{Lemma}
\newtheorem{ppsn}[dfn]{Proposition}
\newtheorem{crlre}[dfn]{Corollary}
\newtheorem{xmpl}[dfn]{Example}
\newtheorem{rmrk}[dfn]{Remark}

\newcommand{\bdfn}{\begin{dfn}\rm}
\newcommand{\bthm}{\begin{thm}}
\newcommand{\blmma}{\begin{lmma}}
\newcommand{\bppsn}{\begin{ppsn}}
\newcommand{\bcrlre}{\begin{crlre}}
\newcommand{\bxmpl}{\begin{xmpl}}
\newcommand{\brmrk}{\begin{rmrk}\rm}

\newcommand{\edfn}{\end{dfn}}
\newcommand{\ethm}{\end{thm}}
\newcommand{\elmma}{\end{lmma}}
\newcommand{\eppsn}{\end{ppsn}}
\newcommand{\ecrlre}{\end{crlre}}
\newcommand{\exmpl}{\end{xmpl}}
\newcommand{\ermrk}{\end{rmrk}}

\newcommand{\bbc}{\mathbb{C}}
\newcommand{\bbz}{\mathbb{Z}}
\newcommand{\bbn}{\mathbb{N}}

\newcommand{\cla}{\mathcal{A}}

\newcommand{\prf}{\noindent{\it Proof\/}: }

\def \qed { \mbox{}\hfill
$\Box$\vspace{1ex}}

\begin{document}

%%%%%%%%%%%%%%%%%%%%%%%%%%%%%%%%%

 \author{\sc  {Partha Sarathi Chakraborty, Bipul Saurabh}}
 \title{Gelfand-Kirillov dimension of the algebra of regular functions on  quantum groups}
 \maketitle

%%%%%%%%%%%%%%%%%%%%%%%%%%%%%%%%%%
%%%%%  ABSTRACT
%%%%%%%%%%%%%%%%%%%%%%%%%%%%%%%%%%
 \begin{abstract} 
 Let  $G_q$ be the $q$-deformation of a simply connected simple compact Lie group $G$ of type $A$,  $C$ or $D$ and $\mathcal{O}_q(G)$ be the  algebra of regular functions on $G_q$. 
In this article,   we prove that the Gelfand-Kirillov dimension  of $\mathcal{O}_q(G)$ is 
 equal to the dimension of real manifold $G$.
\end{abstract}

{\bf AMS Subject Classification No.:} {\large 16}P{\large 90}, {\large
17}B{\large 37}, {\large
  20}G{\large 42}

\newsection{Introduction}   
Motivated by the isomorphism theorem of Weyl algebras, Gelfand and Kirillov \cite{KraLen-2000aa} introduced a 
measure namely, Gelfand Kirillov dimension (abbreviated as GKdim), of growth  of an algebra. 
%In commutative case, it matches with the classical notions. 
For finitely generated
commutative algebra $\cla$, the Gelfand Kirillov dimension is same as the Krull 
dimension of $\cla$ and for commutative domains, it equals the transcedence degree of its fraction 
field (see chapter 4, \cite{KraLen-2000aa}). In order to give a precise estimate for the growth exponent, 
Banica and Vergnioux \cite{BanVer-2009aa} 
proved that for a connected simply connected compact real Lie group $G$,   GKdim of the Hopf-algebra $\mathcal{O}(G)$
generated by matrix co-efficients of all finite dimensional unitary  representations of $G$  is   same as manifold dimension of $G$. In the same article, they mentioned that they do 
not have any other example of Hopf algebra having polynomial growth. 
%Another measure ``$\mbox{dim}$'' introduced by Vergnioux \cite{Ver-2007aa} is related to growth of the vector dimension function on the representation ring $R(G)$ associated with a compact quantum group $G$. 
Later D{\'A}ndrea, Pinzari and Rossi (\cite{DanPinRos-2016aa}) extended their result to   compact Lie groups 
(see Theorem 3.1, \cite{DanPinRos-2016aa}).  
But apart from these commuatative examples, not much is known about the growth of other Hopf algebras. 
For many noncommutative Hopf algebras, even the 
question whether they have a polynomial growth remain unanswered. Therefore it is worthwhile to investigate 
in this direction. The most natural candidate to investigate is  the 
Hopf algebra of finite dimensional unitary  representations of a compact quantum group. In this paper, we take the case of $q$-deformation of a classical Lie group 
of type $A$, $C$ and $D$ and extend the result of Banica and Vergnioux to the noncommutative Peter-Weyl algebra associated with these compact quantum groups.

%we can not expect the same result for the $q$-deformation of a classical Lie group $G$ because  from Corollary 4.5 (\cite{DanPinRos-2016aa}) 
%it follows that polynomial growth of ``$\mbox{dim}$'' will imply 
%the dense Hopf algebra to be of Kac type which is not the case. 
%This raises the question whether the same result is true for %the Hopf-algebra generated by matrix co-efficients of all finite dimensional unitary  representations of the compact quantum group $G_q$
%where $G_q$ is 
%the  $q$-deformation of a semisimple simply connected compact Lie group. In this article, we answer this question in affirmative for $q$-deformation of a classical Lie group 
%of type $A$, $C$ and $D$.% by proving that 
%Gelfand Kirillov dimension  of this Hopf algebra  is equal to  the 
%dimension of real manifold  $G$.

Let $G$ be a semisimple simply connected compact Lie group of rank $n$ and $\mathfrak{g}$ be the Lie algebra of $G$. The algebra of functions $\mathcal{O}_q(G)$ on its  $q$-deformation $G_q$ can be defined as the  
subalgebra of the dual algebra of  quantized universal enveloping algebra $U_q(\mathfrak{g})$ generated by 
matrix co-efficients of all finite dimensional admissible  representations of  $U_q(\mathfrak{g})$. In this article, we are mainly interested in the computation of GKdim for $\mathcal{O}_q(G)$.
In commutative case i.e. for $q=1$, Banica and Vergnioux \cite{BanVer-2009aa} 
proved that the GKdim of polynomial algebra $\mathcal{O}(G)$
is equal to the dimension of $G$ as a real manifold. Here we show that in
noncommutative case i.e for $0 < q <1$, the canonical Hopf $*$-algebra $\mathcal{O}_q(G)$ has GKdim equal to  the  dimension of  real manifold  $G$ provided  $G$ is of type $A$, $C$ and $D$. 
Moreover, we prove that similar results hold for certain quotient spaces of $G_q$ in type $A$ and $C$. 
It  answers 
the query of Banica and Vergnioux by providing  examples of noncommuatative noncocommutative Hopf algebras having polynomial growth.

%Neshveyev and Tuset \cite{NesTus-2012ab} generalises this result to the quotient algebra $C[G_q/H_q]$-modules where $H_q$ is $q$-deformation 
%of a closed Poisson Lie subgroup of $G$. 
%Since Gelfand Kirillov dimension can be defined for a finitely generated module, 
%The computation of Gelfand Kirillov dimension for these algebras and their simple unitarizable modules 
%is  a matter of interest and this brings the content of this article.
%If $G$ is simple simply~connected compact Lie group of type $A$, $B$, $C$ and $D$, then it is well known fact that $\mathcal{O}_q(G)$ is a 
%finitely generated unital algebra. So, to compute Gelfand  Kirillov dimension, one can stick 
%to a fixed generating set.

Let us give a sketch of the proof. We will assume $G$ to be of type $A_n$, $C_n$ or $D_n$. 
Let $\ell(\omega_n)$ be the length of the longest element $\omega_n$ of the Weyl group $W_n$ of $G$.  Let $\mathcal{P}_q(\scrt)$ be 
the algebra of endomorphisms on $c_{00}(\bbn)$ generated by the endomorphisms $\sqrt{1-q^{2N+2}}S$, $\sqrt{1-q^{2N}}S^*$ and $\beta:=q^N$ where $S$ is the left shift operator,  $S^*$ 
is the right shift operator
and $N$ is the number operator. Further let $\mathcal{P}(C(\bbbt))$ be the algebra of endomorphisms $c_{00}(\bbz))$ generated by left shift operator  $S$ and right shift operator $S^*$.  
It can be shown  that the  algebra $\mathcal{O}_q(G)$ can be embedded as a subalgebra of $\mathcal{P}(C(\bbbt))^{\otimes n} \otimes \mathcal{P}_q(\scrt)^{\otimes \ell(\omega_n)}$. 
We show that the algebra generated by the left shift and the right shift in $c_{00}(\bbn)$ has GKdim $2$ and the
algebra generated by the left shift and the right shift in $c_{00}(\bbz)$ has GKdim $1$.  
It proves that GKdim of $\mathcal{O}_q(G)$ is less than $2\ell(\omega_n)+n$. Next, we write $\omega_n$ as a product of $n$ elements in a certain manner and then using recursion 
 we produce  enough number of linealy independent 
endomorphisms in $\mathcal{O}_q(G)$ to get  GKdim of $\mathcal{O}_q(G)$ to be equal to $2\ell(\omega_n)+n$. This completes the proof as the dimension of $G$ as a real manifold is same as $2\ell(\omega_n)+n$.

Organisation of this paper is as follows. Next section is dedicated to the preliminaries on  representation theory of the Hopf $*$-algebra $\mathcal{O}_q(G)$. 
In the third section, we compute GKdim  of the algebra $\mathcal{O}_q(G)$ and prove our main result.  In the final section, we prove similar results for Peter Weyl algebra of some quotient spaces.

%All algebras considered in this article are finitely generated unital algebra over $\bbc$.
%All $\mathcal{O}_q(G)$-modules discussed in the paper are simple unitarizable left $\mathcal{O}_q(G)$-modules.
%\textbf{Notations and conventions:}  
Throughout the paper algebras are assumed to be unital and  over the field $\bbc$.
%Given an algebra $\cla$ and a $\cla$-module $M$, we write GKdim$(\cla)$ and GKdim$(M)$  for  the Gelfand kirillov dimensions of $\cla$ and $M$ respectively. 
Elements of the Weyl group 
will be called Weyl words. We denote by $\ell(w)$  
the length of the Weyl word $w$. 
We used $SP(2n)$ instead of more commonly used notation $SP(n)$ for symplectic group of rank $n$  and hence quantum symplectic group is denoted by $SP_q(2n)$.    Let us denote 
by $\left\{e_n: n\in \bbn \right\}$ and  $\left\{e_n: n\in \bbz \right\}$ the standard bases of the 
vector spaces $c_{00}(\bbn)$ and  $c_{00}(\bbz)$ respectively.
The map $e_n \mapsto e_{n-1}$ will be   denoted by  $S$ and the map  map $e_n \mapsto e_{n+1}$ will be   denoted by  $S^*$.  %For $m< 0$, $(S^*)^m$ denotes the operator $S^{-m}$. 
The map $e_n \mapsto ne_n$ will be called the
number operator  $N$.  We denote by $\overleftarrow\prod_{i=1}^na_i$ the element $a_na_{n-1}\cdots a_1$. % The letter $\scrt$ will denote the Toeplitz algebra. 
%Just to avoid unnecessary endomorphisms, we adopt the following convention.  
Let $T$ and $T^{'}$ be two endomorphisms of the vector space  $c_{00}(\bbz)^{\otimes l}\otimes c_{00}(\bbn)^{\otimes k}$ and $V$ be a subspace of 
$c_{00}(\bbz)^{\otimes l}\otimes c_{00}(\bbn)^{\otimes k}$.  We say that 
$T \sim T^{'}$ on $V$ if there exist natural numbers $m_1,m_2, \cdots ,m_{k}$ and a nonzero constant $C$ such that 
\[
 T=CT^{'}(\underbrace{1 \otimes 1 \otimes \cdots \otimes 1}_{l \mbox{ copies}} \otimes q^{m_1N}\otimes q^{m_2N}\otimes \cdots \otimes q^{m_{k}N})
\]
on $V$. Throughout this paper, $q$ will denote a real number in the interval  $(0,1)$ and $C$ is used to denote a generic constant.

\newsection{Quantized  algebra of regular functions}
In this section, we recall the definition of quantized algebra of regular functions on a simply connected semisimple compact Lie group $G$ and give a faithful homomorphism of this algebra in order to find a new set of generators 
consisting of endomorphisms of a vector space. 
For a detailed treatment, we refer the reader to (\cite{KliSch-1997aa}, Chapter $3$ in \cite{KorSoi-1998aa}). 
Let $G$ be a simply connected semisimple compact Lie group of rank $n$ and  $\mathfrak{g}$ be its complexified Lie algebra. Fix a nondegenerate symmetric ad-invariant form  $\left\langle\cdot,\cdot\right\rangle$ on $\mathfrak{g}$ 
such that its restriction to the real Lie algebra of $G$ is negative definite.
Let $\Pi:=\{ \alpha_1,\alpha_2, \cdots, \alpha_n \}$   be the set of simple roots. % Here indices of the roots are according to the standard order. 
For simplicity, we write the root $\alpha_i$ as $i$ and the reflection $s_{\alpha_i}$ defined by the root $\alpha_i$ as $s_i$. % for all $1 \leq i \leq n$. 
The  Weyl group $W_n$ of $G$ can be described as the group generated 
 by the reflections $\{s_{i}: 1\leq i \leq n\}$. 
\bdfn  \label{chap2-d-G_q}
Let $U_q(\mathfrak{g})$ be the quantized universal enveloping algebra of $\mathfrak{g}$.   It has a 
%compact real form 
$*$-structure corresponding to the compact real form of $\mathfrak{g}$ (see page $161$ and $179$, \cite{KliSch-1997aa}).  
%Then the quantized algebra of regular functions $\mathcal{O}_q(G)$ 
The Hopf $*$-subalgebra of the dual Hopf $*$-algebra of  $U_q(\mathfrak{g})$ consisting of matrix co-efficients
of finite dimensional unitarizable $U_q(\mathfrak{g})$-modules is called 
the quantized algebra of regular functions on $G$ (see page $96-97$, \cite{KorSoi-1998aa}).   It is denoted by $\mathcal{O}_q(G)$. %Let $C(G_q)$ be the enveloping $C^*$-algebra
 %of the Hopf $*$-algebra $\mathcal{O}_q(G)$. Then extending the comultiplication map $\Delta$ to $C(G_q)$,  we get a 
%compact quantum group $(C(G_q),\Delta)$ which is called  the
%$q$-deformation of the  compact  Lie group $G$. It is denoted by $G_q$.
\edfn

Let $(\!(u_{j,\mathfrak{g}}^i)\!)$ be the defining corepresentation of $\mathcal{O}_q(G)$ if $G$ is of type $A_n$ and $C_n$ and  
the irreducible corepresentation of $\mathcal{O}_q(G)$ corresponding to the highest weight $(1,0,0,\cdots ,0)$ if $G$ is of type  $D_n$. 
In first case, entries of the matrix $(\!(u_{j,\mathfrak{g}}^i)\!)$ generate the Hopf $*$-algebra $\mathcal{O}_q(G)$. %  defining representation of $G_q$ if $G_q=SU_q(n+1)$ or $SP_q(2n)$ and  
%the irreducible representation of $G_q$ corresponding to the highest weight $(1,0,0,\cdots ,0)$ if $G_q=\mbox{Spin}_q(n)$. We denote the dimension of the representation $(\!(u_j^i)\!)$ by $N_n$.
 %For $SU_q(n+1)$ or $SP_q(2n)$, the algebra $\mathcal{O}_q(G)$ is generated by $u_j^i$'s.
In latter case, they genarate a proper Hopf $*$-subalgebra of $\mathcal{O}_q(\mbox{Spin}(2n))$ which we denote as $\mathcal{O}_q(SO(2n))$.  
% (see page 410, Theorem 22, \cite{KliSch-1997aa}). 
The generators of  $\mathcal{O}_q(\mbox{Spin}(2n))$ are the  matrix entries of the corepresentation 
$(\!(z_j^i)\!)$ of $\mathcal{O}_q(\mbox{Spin}(2n))$  with highest weight $(1/2,1/2,\cdots ,1/2)$. 
We denote the dimension of the corepresentation $(\!(u_{j,\mathfrak{g}}^i)\!)$ by $N_n$.  We will 
drop the subscript $\mathfrak{g}$ in $(\!(u_{j,\mathfrak{g}}^i)\!)$ whenever the Lie algebra 
$\mathfrak{g}$ is clear from the context.
Using a result of Korogodski and Soibelman (\cite{KorSoi-1998aa}), we will now describe  all simple unitarizable $\mathcal{O}_q(G)$-modules.

\textbf{Elementary simple unitarizable $\mathcal{O}_q(G)$-modules:} Let $d_i=\left\langle \alpha_i, \alpha_i \right\rangle /2$ and $q_i=q^{d_i}$ for $1 \leq i \leq n$.  
Define $\phi_i :U_{q_i}(\mathfrak{sl}(2)) \longrightarrow  U_q(\mathfrak{g})$
be a  $*$-homomorphism given on the generators of $U_{q_i}(\mathfrak{sl}(2))$ by, 
\begin{displaymath}
 K \longmapsto K_i, \qquad E \longmapsto E_i, \qquad F \longmapsto F_i.
\end{displaymath}
By duality, it induces an  epimorphism 
\begin{displaymath}
\phi_i^* : \mathcal{O}_q(G) \longrightarrow \mathcal{O}_{q_i}(SU(2)).
\end{displaymath}
We  will use this map to get all elementary simple unitarizable modules of $\mathcal{O}_q(G)$. % Let $N$ be the number operator given by $N: e_{n} \mapsto n e_{n}$ 
%and $S$ be the  shift operator given by $S: e_{n} \mapsto e_{n-1}$ on  $c_{00}(\bbn)$. 
Denote by $\Psi$ the following 
action of $\mathcal{O}_q(SU(2))$ on $c_{00}(\bbn)$ (see Proposition $4.1.1$, \cite{KorSoi-1998aa});
\begin{IEEEeqnarray}{rCl} \label{su(2)}
\Psi(u_l^k)e_p=\begin{cases}
              \sqrt{1-q^{2p}}e_{p-1} & \mbox{ if } k=l=1,\cr
              \sqrt{1-q^{2p+2}}e_{p+1} & \mbox{ if } k=l=2,\cr
							-q^{p+1}e_p & \mbox{ if } k=1,l=2,\cr
							q^pe_p & \mbox{ if } k=2,l=1,\cr
							\delta_{kl}e_p & \mbox{ otherwise }.\cr
              \end{cases}
\end{IEEEeqnarray}
For each $1 \leq i \leq n$, define an action  $\pi_{s_{i}}^n := \Psi \circ \phi_{i}^{*}$ of $\mathcal{O}_q(G)$. %Each $\pi_{s_{i}}^n$ is an  irreducible representation and is called an elementary representation of $\mathcal{O}_q(G)$. 
Each $\pi_{s_{i}}^n$ gives rise to 
an elementary simple $\mathcal{O}_q(G)$-module $V_{s_i}$. Also, for each $t \in \bbbt^n$, there are one dimensional  $\mathcal{O}_q(G)$-module $V_t$ with the action $\{\tau_t^n\}$.  
%The simple $\mathcal{O}_q(G)$-modules associated with the representation $\tau_t^n$ is denoted by $V_t$.
Given two actions  $\varphi$ and $\psi$ of $\mathcal{O}_q(G)$, define an action $\varphi * \psi := (\varphi \otimes \psi)\circ \Delta$. Similarly for any two $\mathcal{O}_q(G)$-module $V_{\varphi}$ and $V_{\phi}$, define 
$V_{\varphi} \otimes V_{\phi}$ as $\mathcal{O}_q(G)$-module with $\mathcal{O}_q(G)$ action coming from  $\varphi * \psi$.
For $w \in W_n$ such that  $s_{i_{1}}s_{i_{2}}...s_{i_{k}}$ is a reduced expression for $w$ and $t \in \bbbt^n$, 
%define an action
%$\pi_{w}^n$ as  $\pi_{s_{i_{1}}}^n*\pi_{s_{i_{2}}}^n*\cdots *\pi_{s_{i_{k}}}^n$  and  the corresponding $\mathcal{O}_q(G)$-module as $V_{w}$.  
%Then $V_{w}$ is an simple unitarizable which is independent of the reduced expression. %Besides these representations, there are one dimensional irreducible representations $\{\tau_t^n\}_{t \in \bbbt^n}$ parametrized by $\bbbt^n$.  
 %Moreover, for $t \in \bbbt^{n}, w \in W$,  
 define an action $\pi_{t,w}^n$ by $\tau_{t}^n*\pi_{s_{i_{1}}}^n*\pi_{s_{i_{2}}}^n*\cdots *\pi_{s_{i_{k}}}^n$ and denote the corresponding $\mathcal{O}_q(G)$-module by $V_{t,w}$. If $t=1$, we write the 
 action $\pi_{t,w}^n$ as $\pi_w^n$ and the associated module  $V_{t,w}$ by $V_w$. 
 We refer the reader to (\cite{KorSoi-1998aa}, page~121) for the following theorem.
\bthm \label{allsimple}
The set $\left\{V_{t,w}; t \in \bbbt^{n}, w \in W_n\right\}$ is a complete set of mutually inequivalent simple  unitarizable left $\mathcal{O}_q(G)$-module.
\ethm
Define the endomorphisms $\alpha:=\sqrt{1-q^{2N+2}}S$, $\alpha^*:=\sqrt{1-q^{2N}}S^*$ and $\beta:=q^N$ acting on the vector space $c_{00}(\bbn)$.
Let $\mathcal{P}_q(\scrt) \subset$ END$(c_{00}(\bbn))$ be the algebra generated by $\alpha$, $\alpha^*$ and $\beta$ and 
$\mathcal{P}(C(\bbbt))\subset$ END$(c_{00}(\bbz))$ be the algebra generated by $S$ and $S^*$. Given  a Weyl word $w$ of length $\ell(w)$,  we define a homomorphism 
$\chi_{w}^n :\mathcal{O}_q(G)\longrightarrow \mathcal{P}(C(\bbbt))^{\otimes n}\otimes \mathcal{P}_q(\scrt)^{\otimes \ell(w)}$ 
such that $\chi_{w}^n(a)(t) = \pi_{t,w}^n(a)$  for all $a\in \mathcal{O}_q(G)$. 
%By mapping $z \mapsto z^m$ with $e_m$, one can identify $L_2(\bbbt)$ with $L_2(\bbz)$ and hence one 
%can view $C(\bbbt)$ as a subset of $\mathcal{L}(L_2(\bbz))$ by identifying  generator $z \mapsto z$ of $C(\bbbt)$ acting as a  multiplication operator on  $L_2(\bbbt)$ with  the operator $S^*$ on $L_2(\bbz)$. 
%Hence we have  $C(\bbbt^{n})\otimes \scrt^{\otimes \ell(w)} \subset \mathcal{L}(L_2(\bbz)^{\otimes n} \otimes L_2(\bbn)^{\otimes \ell(w)})$ and therefore we will 
%view $\chi_{w}$ as a representation of $C[G_q]$ on the Hilbert space $L_2(\bbz)^{\otimes n} \otimes L_2(\bbn)^{\otimes \ell(w)}$. The following theorem gives a faithful representation of 
%$C[G_q]$.
\bthm \label{faithful}
Let $\omega_{n}$ be the longest word of the Weyl group of $G$. 
Then the homomorphism
\[
                            \chi_{\omega_{n}}^n : \mathcal{O}_q(G)\longrightarrow \mathcal{P}(C(\bbbt))^{\otimes n}\otimes \mathcal{P}_q(\scrt)^{\otimes \ell(\omega_n)}
\]
is faithful.
\ethm
\prf Consider the  enveloping $C^*$-algebra $C(G_q)$ of the Hopf $*$-algebra  $\mathcal{O}_q(G)$. For each $w \in W_n$ and $t \in \bbbt^n$, one can extend the irreducible representation $\pi_{t,w}^n$ and homomorphism $\chi_w^n$ 
to the $C^*$-algebra $C(G_q)$ which we will denote by the same symbols. It follows from \cite{NesTus-2012ab} that the set  $\left\{\pi_{t,w}^n; t \in \bbbt^{n}, w \in W_n\right\}$
is a complete set of mutually inequivalent irreducible representations of $C(G_q)$. It is not difficult to show that if $w^{'}$ is a subword of $w$ then the representation $\pi_{t,w^{'}}^n$ factors through 
the homomorphism $\chi_w^n$. Since  $\omega_{n}$ is the longest word of $W_n$, it follows that each irreducible representation factors through  $\chi_{\omega_{n}}^n$. As a consequence, 
the homomorphism $\chi_{\omega_{n}}^n: C(G_q) \rightarrow C(\bbbt^n)\otimes \scrt^{\otimes \ell(\omega_n)}$ is faithful. Restricting this homomorphism to the subalgebra $\mathcal{O}_q(G)$ 
proves the claim.
\qed\\
Consider the action $\chi_e^n$ of $\mathcal{O}_q(G)$ on the vector space $c_{00}(\bbz)^{\otimes n}$. It is not difficult to see that $\chi_e^n(a)(t)=\tau_t(a)$ for all $a \in \mathcal{O}_q(G)$. 
Therefore for any $w \in W$, we have  $\chi_{w}^n=\chi_e^n * \pi_{w}^n$. We will  explicitly write down the endomorphisms $\chi_e^n(u_j^i)$ of $\mathcal{O}_q(G)$ for type $A_n$, $C_n$ or $D_n$.\\
For $\mathcal{O}_q(G)=\mathcal{O}_q(SU(n+1))$, 
\begin{IEEEeqnarray}{rCl}
 \chi_e^n(u_j^i)=\begin{cases}
          \delta_{ij}1 \otimes 1 \otimes \cdots 1 \otimes \underbrace{S^*}_{n+2-i^{\mbox{th place}}} \otimes 1\otimes \cdots \otimes 1 & \mbox{ if } i \neq 1, \cr
          \delta_{ij}S \otimes S \otimes \cdots \otimes S &\mbox{ if } i =1.   \cr       
        \end{cases} \nonumber
\end{IEEEeqnarray}
For $\mathcal{O}_q(G)=\mathcal{O}_q(SP(2n))$ or $\mathcal{O}_q(\mbox{Spin}(2n))$, 
\begin{IEEEeqnarray}{rCl}
 \chi_e^n(u_j^i)=\begin{cases}
          \delta_{ij}1 \otimes 1 \otimes \cdots 1 \otimes \underbrace{S^*}_{2n+1-i^{\mbox{th place}}} \otimes 1\otimes \cdots \otimes 1 & \mbox{ if } i > n, \cr
          \delta_{ij}1 \otimes 1 \otimes \cdots 1 \otimes \underbrace{S}_{i^{\mbox{th place}}} \otimes 1\otimes \cdots \otimes 1 &\mbox{ if } i \leq n.   \cr       
        \end{cases} \nonumber
\end{IEEEeqnarray}
Looking at the expression of $\chi_e^n(u_j^i)$, it follows that 
\begin{IEEEeqnarray}{rCl} \label{Gamma}
  \chi_w^n(u_j^i)=(\chi_e^n \otimes \pi_w^n)(\Delta(u_j^i))=(\chi_e^n\otimes \pi_w^n)(\sum_{k=1}^{N_n}u_k^i\otimes u_j^k)=\chi_e^n(u_i^i)\otimes \pi_w^n(u_j^i).
\end{IEEEeqnarray}

%%%%%%%%%%%%%%%%%%%%%%%%%%%%%%%%%%%%%%%%%%%%%%%

\newsection{Main result}
In the present section, we show that Gelfand-Kirillov dimension of quantized algebra of regular functions on a simply connected simple compact Lie group $G$ of type $A$,  $C$ or $D$ is 
equal to the dimension of $G$ as a real manifold. Unless otherwise specified,  we denote by  
$\mathcal{O}_q(G)$  one of the Hopf $*$-algebras $\mathcal{O}_q(SU(n+1))$, $\mathcal{O}_q(SP(2n))$ or $\mathcal{O}_q(\mbox{Spin}(2n))$.
\bdfn (\cite{KraLen-2000aa})
Let $A$ be a  unital algebra. The Gelfand-Kirillov dimension of $A$ is given by 
\[
 \mbox{GKdim}(A)=sup_{V}\varlimsup\frac{\ln \dim(V^k)}{\ln k}  
\]
where the supremum is taken over all finite dimensional subspace  $V$ of $A$ containing $1$.  If 
 $A$ is a finitely generated unital algebra then
\[
 \mbox{GKdim}(A)=sup_{\xi}\varlimsup\frac{\ln \dim(\xi^k)}{\ln k}  
\]
where the supremum is taken over all finite sets  $\xi$  containing $1$ that generates $A$. 
\edfn
\brmrk
The quantity ``$\varlimsup\frac{\ln \dim(\xi^k)}{\ln k}$'' does not depend on particular choices of $\xi$  and hence one can choose a fixed (but finite) set of generators of $A$.
 % \item Since we will be dealing with finitely generated algebras, we 
 %defined Gelfand kirillov dimension for this case only. For general case, we refer the reader to (\cite{DanPinRos-2016aa}, \cite{KraLen-2000aa}). 
\ermrk
We state some properties of Gelfand-Kirillov dimension omitting their straightforward proofs.
\begin{itemize}
 \item If $B$ is a finitely generated unital subalgebra of  $A$ then $GKdim(B)\leq GKdim(A)$ (see \cite{Row-1990aa}).
 \item %If $A$ and $B$ are unital then 
 $GKdim(A\otimes B)\leq GKdim(A)+GKdim(B)$.
 %\item $GKdim(A[x])=GKdim(A)+1$.
 %\item If $A$ is finite dimensional algebra then $GKdim(A)=0$.
 %\item If $M$ is finitely generated left $A$-module then $GKdim(M)\leq GKdim(A)$.
\end{itemize}
%Define the endomorphisms $\alpha:=\sqrt{1-q^{2N+2}}S$ and $\beta:=q^N$ acting on the vector space $c_{00}(\bbn)$. Let $\mathcal{P}_q(\scrt) \subset \mathcal{L}(c_{00}(\bbn))$ be the algebra generated by $\alpha$, $\alpha^*$ and $\beta$ and 
%$\mathcal{P}(C(\bbbt))\subset \mathcal{L}(c_{00}(\bbz))$ be the algebra generated by $S$ and $S^*$. Then we have
\bppsn \label{gktau}
GKdim$(\mathcal{P}(C(\bbbt)))=1$ and GKdim$(\mathcal{P}_q(\scrt) )=2$.
\eppsn
\prf Clearly $\{1,S,S^*\}^m=\mbox{span}\{S^k:-m\leq k \leq m\}$ and hence GKdim$(\mathcal{P}(C(\bbbt)))=1$. To show the other claim, take the generating set of $\mathcal{P}_q(\scrt)$ to be $F=\{1,\alpha,\alpha^*,\beta \}$. 
From the commutation relations $q\beta \alpha=\alpha \beta$ and $\alpha \alpha^* -\alpha^* \alpha =(1-q^2)\beta^2$, it is easy to see that 
\[
 F^m=\mbox{span}\{(\alpha^*)^{m_1}\beta^{m_2}\alpha^{m_3}:m_1+m_2+m_3\leq m\}.
\]
Since $\beta^2=1-\alpha^*\alpha$, we get 
\[
 F^m=\mbox{span}\Big\{\{(\alpha^*)^{m_1}\beta\alpha^{m_3}:m_1+m_3 < m\} \cup \{(\alpha^*)^{m_1}\alpha^{m_3}:m_1+m_3 \leq m\} \Big\}.
\]
Hence the dimension of $F^m$ is less than or equal to   $(m+1)^2$. %Morever, 
%\[
 % F^m \supset \mbox{span}\{(\alpha^*)^{m_1}\alpha^{m_3}:m_1+m_3\leq m\}.
%\]
Since $\{(\alpha^{*})^{m_1}\alpha^{m_3}:m_1+m_3 \leq  m\}$ are linearly independent set of endomorphisms, we conclude that the  dimension of $F^m$ is greater than or equal to   $m+1 \choose 2$. 
Putting together, we get  GKdim$(\mathcal{P}_q(\scrt) )=2$.
\qed
\blmma \label{lessthandim} 
Let $\omega_n$ be  the longest element of the Weyl group  of $G$. Then  one has 
\[
 \mbox{GKdim}(\mathcal{O}_q(G))\leq 2\ell(\omega_n)+n.
\]
\elmma
\prf
By Theorem \ref{faithful}, the algebra  $\mathcal{O}_q(G)$ can be viewed as a subalgebra of $\mathcal{P}(C(\bbbt))^{\otimes n} \otimes \mathcal{P}_q(\scrt)^{\otimes \ell(\omega_n)}$. Using
the properties of Gelfand-Kirillov dimension mentioned above and  Proposition \ref{gktau}, we have
\begin{IEEEeqnarray}{rCl}
 GKdim(\mathcal{O}_q(G))\leq GKdim(\mathcal{P}(C(\bbbt))^{\otimes n} \otimes \mathcal{P}_q(\scrt)^{\otimes \ell(\omega_n)})\leq 2\ell(\omega_n)+n. \nonumber
\end{IEEEeqnarray}
This settles the claim.
\qed\\
In what follows, we will show that equality holds in Lemma \ref{lessthandim}. Our strategy is similar to that  given in \cite{ChaSau-2017aa} with some modifications.
First we need the following result.
\blmma \label{independent}
Let  $\alpha_d$ be the endomorphism $S\sqrt{1-q^{2dN}}$ of the vector space $c_{00}(\bbn)$. For any  fixed $j, k \in \bbn$ and $ 0 \leq i \leq j$, let
$T_i \sim  \alpha_d^i(\alpha_d^*)^{i+k}$ on $c_{00}(\bbn)$. Then elements of the set $\{T_i: 0 \leq i \leq j\}$ 
are linearly independent endomorphisms.
\elmma
\prf 
Enough to prove for $k=0$. Consider the set
 $\{q^{a_rN}:  1 \leq r \leq s, \mbox{ no two } a_r\mbox{'s are same}\}$.  Let $V$ be the following Vandermonde matrix; 
\[V= 
  \left[ {\begin{matrix}
   1 & q^{a_1} & q^{2a_1}& \cdots & q^{(s-1)a_1}\\
    1 & q^{a_2} & q^{2a_2}& \cdots & q^{(s-1)a_2}\\
    \cdot & \cdot & \cdot & \cdots &\cdot \\
    \cdot & \cdot & \cdot & \cdots &\cdot \\
     1 & q^{a_s} & q^{2a_s}& \cdots & q^{(s-1)a_s}\\
  \end{matrix} } \right].
\]
Since $\det(V)=\prod_{p\neq r}(q^{a_p}-q^{a_r})\neq 0$, it follows that elements of 
the set $\{q^{a_rN}:  1 \leq r \leq s, \mbox{ no two } a_r\mbox{'s are same}\}$ are linealy independent. 
Next, we have
\begin{IEEEeqnarray}{rrCl}
                   &T_i &\sim& \alpha_d^i(\alpha_d^*)^{i} 
                   \sim (1-q^{2dN+2d}) (1-q^{2dN+4d}) \cdots (1-q^{2dN+2dm}), \nonumber \\
                   \Rightarrow&T_i&= &C_iq^{m_iN}(1-q^{2dN+2d}) \cdots (1-q^{2dN+2dm}) \nonumber 
 \end{IEEEeqnarray}
for some $m_i \in \bbn$ and nonzero constant $C_i$. If we expand the right hand side, we get $2^m$ terms of the form $q^{b_r}q^{a_rN}$ such that all $a_r$'s are different. Let us assume that  
$\sum_{i=1}^{j}c_iT_i=0$. Since $T_j$ has $2^j$  terms of the form $q^{a_rN}$ and $2^j > \sum_{i=1}^{j-1}2^i$, we get $c_jq^{sN}=0$ for some $s \in \bbn$ which further implies 
that $c_j=0$. Repeating the same argument, we get $c_i=0$ for all $1\leq i \leq j$ and this completes the proof.
\qed \\
% For that, we need to extend Lemma \ref{poly} to the representations of type $\chi_{w}$ for $w \in W$. 
%For $w \in W_n$ and  $\ell_i:=\ell(w_i)$.
We will recall from \cite{ChaSau-2017aa} some results that will be needed to prove our main claim. 
\blmma \label{sample1}
% Let $\mathcal{O}_q(G)$ be one of the Hopf $*$-algebras $\mathcal{O}_q(SU(n+1))$, $\mathcal{O}_q(SP(2n))$ or $\mathcal{O}_q(\mbox{Spin}(2n))$.
Let $w \in W_n$. 
Then there exist polynomials $p_{1}^{(w,n)}, p_{2}^{(w,n)},\cdots p_{\ell(w_n)}^{(w,n)}$ with noncommuting variables $\pi_w^n(u_j^{N_n})$'s 
and a permutation $\sigma$ of $\{1,2,\cdots ,\ell(w_n)\}$ such that  for all $r_1^n,r_2^n,\cdots,r_{\ell(w_n)}^n \in \bbn$, one has
 \begin{IEEEeqnarray}{lCl} \label{so2n}
 p_{j}^{(w,n)} \sim 1^{\otimes \sum_{i=1}^{n-1} \ell(w_i)}\otimes 1^{\otimes \sigma(j)-1}  \otimes \sqrt{1-q^{2d_{\sigma(j)}N}}S^*\otimes 1^{\otimes \ell(w_n)- \sigma(j)} 
\end{IEEEeqnarray}
on the subspace generated by standard  basis elements having $e_0$ at $(\sum_{i=1}^{n-1}\ell(w_i)+\sigma(k))^{th}$ place for $k <j$.
\elmma
\prf
See the proof of the Lemma $3.4$ in \cite{ChaSau-2017aa}.
\qed \\
An element $w$ of $W_n$ can be written  in a reduced form as: 
$\psi_{1,k_{1}}^{(\epsilon_{1})}(w)\psi_{2,k_{2}}^{(\epsilon_{2})}(w)\cdots \psi_{n,k_{n}}^{(\epsilon_{n})}(w)$ for some choices 
of $\epsilon_1, \epsilon_2, \cdots, \epsilon_n$ and $k_1,k_2,\cdots k_n$ where
$\epsilon_{r} \in \left\{0,1,2\right\}$ and $n-r+1\leq k_{r}\leq n$ with the convention that, \\
\textbf{Case 1}: $\mathfrak{sl}(n+1)$
\[
 \psi_{r,k_{r}}^{\epsilon}(w)= \begin{cases}
               s_{r}s_{r-1}\cdots s_{n-k_r+1} & \mbox{ if } \epsilon=1,2\cr
							 \mbox{ empty string } & \mbox{ if } \epsilon=0.\cr
               \end{cases}
\]
\textbf{Case 2}: $\mathfrak{sp}(2n)$
\[
\psi_{r,k_{r}}^{\epsilon}(w)=\begin{cases}
               s_{n-r+1}s_{n-r+2}\cdots s_{k_r} & \mbox{ if } \epsilon=1,\cr
               s_{n-r+1}s_{n-r+2}\cdots...s_{n-1}s_{n}s_{n-1}\cdots s_{k_{r}} & \mbox{ if } \epsilon=2,\cr
							 \mbox{ empty string } & \mbox{ if } \epsilon=0.\cr
               \end{cases}
\]
\textbf{Case 3}: $\mathfrak{so}(2n)$
\[
\psi_{r,k_{r}}^{\epsilon}(w)=\begin{cases}
               s_{n-r+1}s_{n-r+2}\cdots s_{k_r} & \mbox{ if } \epsilon=1,\cr
               s_{n-r+1}s_{n-r+2}\cdots...s_{n-1}s_{n}s_{n-2}s_{n-3}\cdots s_{k_{r}} & \mbox{ if } \epsilon=2,\cr
							 \mbox{ empty string } & \mbox{ if } \epsilon=0.\cr
               \end{cases}
\]
For details, we refer the reader to (\cite{Hum-1990aa} or subsection 2.2 in \cite{ChaSau-2017aa}). 
We  call the word $\psi_{r,k_{r}}^{(\epsilon_{r})}(w)$ 
the $r^{th}$ part $w_r$ of $w$. 
For type $C_n$ and $D_n$, let  $M_n^i=n-i+1$  and $N_n^i=N_n-n+i$ and for type $A_n$, let
$M_n^i=1$ and $N_n^i=i+1$.
\blmma \label{unique path}% Let $\mathcal{O}_q(G)$ be one of the Hopf $*$-algebras $\mathcal{O}_q(SU(n+1))$, $\mathcal{O}_q(SP(2n))$ or $\mathcal{O}_q(\mbox{Spin}(2n))$.
  Let $w \in W_n$ be of the form   $w=w_{i+1}w_{i+2}\cdots w_l$, $l\leq n$ and let  $V_w$ be the associated   $\mathcal{O}_q(G)$-module.
  Then for each $M_n^i \leq k \leq N_n^i$, there exists unique $r_w(k) \in \{M_n^{i+l},\cdots ,N_n^{i+l}\}$ 
  such that 
\[
\pi_w^n(u_{r_w(k)}^k)(e_0\otimes e_0\otimes \cdots \otimes e_0)=Ce_0\otimes e_0\otimes \cdots \otimes e_0
\] 
 where $C$ is a nonzero real number. 
Moreover,
\begin{enumerate}
 \item 
 $\pi_w^n(u_{r_w(k)}^j)(e_0\otimes e_0\otimes \cdots \otimes e_0)=0$ for $j \in \{M_n^i,\cdots ,N_n^i\}/\{k\}$.
 \item
 $\pi_w^n((u_{r_w(k)}^k)^*)(e_0\otimes e_0\otimes \cdots \otimes e_0)=Ce_0\otimes e_0\otimes \cdots \otimes e_0$.
 \item
  $\pi_w^n((u_{r_w(k)}^j)^*)(e_0\otimes e_0\otimes \cdots \otimes e_0)=0$ for $j \in \{M_n^i,\cdots ,N_n^i\}/\{k\}$. 
\end{enumerate}
\elmma
\prf In Lemma $3.6$ in  \cite{ChaSau-2017aa}, explicit description of the function $r_w(k)$ is given. Using that and diagram representation,  one can verify 
the claim.
\qed\\
The following lemma gives some linearly independent endomorphisms   which when applied  to a fixed vector of the form  $e_0 \otimes v \otimes e_0^{\otimes \ell(w_n)}$ 
give all matrix units of the the first component and the components in  the $n^{th}$ part. More precisely, 
\blmma \label{sample3}
%Let $\mathcal{O}_q(G)$ be one of the Hopf $*$-algebras $\mathcal{O}_q(SU(n+1))$, $\mathcal{O}_q(SP(2n))$ or  $\mathcal{O}_q(\mbox{Spin}(2n))$. 
Let $w \in W_n$. Then there exist a permutation $\sigma$ of $\{1,2,\cdots ,\ell(w_n)\}$ and polynomials $g_{0}^{(w,n)}, g_{1}^{(w,n)},\cdots ,g_{\ell(w_n)}^{(w,n)}$ and  $g_{1*}^{(w,n)},\cdots ,g_{\ell(w_n)*}^{(w,n)}$ with
variables $\chi_w^n(u_j^{N_n})$'s and $\chi_w^n((u_j^{N_n})^*)$'s 
such that 
\begin{enumerate}
 \item
 \begin{IEEEeqnarray}{lCl} 
 (g_{0}^{(w,n)})^{r_0^n}(g_{1*}^{(w,n)})^{s_{\sigma(1)}^n}(g_{1}^{(w,n)})^{r_{\sigma(1)}^n}
 (g_{2*}^{(w,n)})^{s_{\sigma(2)}^n}(g_{2}^{(w,n)})^{r_{\sigma(2)}^n}\cdots (g_{\ell(w_n)*}^{(w,n)})^{s_{\sigma(\ell(w_n))}^n}  \nonumber \\ 
 (g_{\ell(w_n)}^{(w,n)})^{r_{\sigma(\ell(w_n))}^n}  (e_0\otimes v \otimes e_0^{\otimes \ell(w_n)}) 
 = C e_{r_0^n}\otimes  v \otimes  e_{r_{1}^n-s_1^n}\otimes \cdots \otimes e_{r_{\ell(w_n)}^n-s_{\ell(w_n)}^n} \nonumber 
\end{IEEEeqnarray}
where  $v \in c_{00}(\bbz)^{\otimes n-1}\otimes  c_{00}(\bbn)^{\otimes \sum_{k=1}^{n-1}\ell(w_k)}$, $r_0^n \in \bbn, r_i^n, s_i^n \in \bbn$ and $r_i^n\geq s_i^n$ for $1 \leq i \leq \ell(w_n)$.
 \item The elements of the set
 \begin{IEEEeqnarray}{rCl}
   \Big\{ (g_{0}^{(w,n)})^{r_0^n}(g_{1*}^{(w,n)})^{s_{\sigma(1)}^n}(g_{1}^{(w,n)})^{r_{\sigma(1)}^n} (g_{2*}^{(w,n)})^{s_{\sigma(2)}^n}(g_{2}^{(w,n)})^{r_{\sigma(2)}^n} \cdots \nonumber \\
 \cdots (g_{\ell(w_n)*}^{(w,n)})^{s_{\sigma(\ell(w_n))}^n} (g_{\ell(w_n)}^{(w,n)})^{r_{\sigma(\ell(w_n))}^n}: r_0^n \in \bbn, r_i^n, s_i^n \in \bbn, r_i^n\geq s_i^n \in \bbn \Big\} \nonumber 
 \end{IEEEeqnarray}
are linearly independent  endomorphisms.  
\end{enumerate}
\elmma
\prf Define 
\[
 g_0^{(w,n)}:=\chi_{w}^n(u_{N_n-\ell(w_n)}^{N_n}).
\]
Since the endomorphism  $\pi_{w}^n(u_{N_n-\ell(w_n)}^{N_n})$ is of the form $1^{\otimes \sum_{k=1}^{n-1}\ell(w_k)} \otimes q^{m_1N}\otimes q^{m_2N}\otimes \cdots \otimes q^{m_{\ell(w_n)}N}$, we get 
\begin{IEEEeqnarray}{lCl} \label{eq1}
 g_0^{(w,n)}&=& \chi_{w}^n(u_{N_n-\ell(w_n)}^{N_n}) 
              =\chi_e^n * \pi_{w}^n(u_{N_n-\ell(w_n)}^{N_n}) 
              = \chi_e^n(u_{N_n}^{N_n}) \otimes \pi_{w}^n(u_{N_n-\ell(w_n)}^{N_n}) \nonumber \\
             % &=& \underbrace{S^* \otimes 1 \otimes \cdots \otimes 1}_{n \mbox{ times}} \otimes q^{d_1N}\otimes q^{d_2N}\otimes \cdots \otimes q^{d_{\ell(w_n)}N} \nonumber 
             &\sim& \underbrace{S^* \otimes 1 \otimes \cdots \otimes 1}_{n \mbox{ times}} \otimes 1^{\otimes \sum_{k=1}^{n}\ell(w_k)}
\end{IEEEeqnarray}
on the whole vector space. 
%for some  $m_i \in \bbn$ for $1 \leq i \leq \ell(w_n)$. 
For $1\leq j \leq \ell(w_n)$, let  $h_j^{(w,n)}$ be the polynomial obtained by replacing the action $\pi_{w}^n$ 
in the polynomial $p_j^{(w,n)}$ given in Lemma \ref{sample1} with $\chi_{w}^n$. Define 
\[
 g_j^{(w,n)} := (g_0^{(w,n)})^*)^{s}h_j^{(w,n)}
\]
 where $s =\mbox{ degree of }p_j^{(w,n)}$.  Hence we have
\begin{IEEEeqnarray}{lCl} \label{eq2}
 g_j^{(w,n)}= (g_0^{(w,n)})^*)^{s}h_j^{(w,n)} \sim \underbrace{1 \otimes 1 \otimes \cdots \otimes 1}_{n \mbox{ times}} \otimes p_j^{(w,n)}
\end{IEEEeqnarray}
on the subspace generated by standard orthonormal basis elements having $e_0$ at $(n+\sum_{i=1}^{n-1}\ell(w_i)+\sigma(k))^{th}$ place for $k <j$. Define $g_{j*}^{(w,n)}=(g_j^{(w,n)})^{*}$.
Now part (1) of the claim follows from  Lemma \ref{sample1}. Further it follows from first part of the claim that possible dependency can occur among the elements of the type:
\begin{IEEEeqnarray}{rCl}
   \{ (g_{0}^{(w,n)})^{r_0^n}(g_{1*}^{(w,n)})^{p_{\sigma(1)}^n}(g_{1}^{(w,n)})^{r_{\sigma(1)}^n} (g_{2*}^{(w,n)})^{p_{\sigma(2)}^n}(g_{2}^{(w,n)})^{r_{\sigma(2)}^n} \cdots \nonumber \\
 \cdots (g_{\ell(w_n)*}^{(w,n)})^{p_{\sigma(\ell(w_n))}^n} (g_{\ell(w_n)}^{(w,n)})^{r_{\sigma(\ell(w_n))}^n}:  r_i^n - p_i^n= c_i^n  \}. \nonumber 
 \end{IEEEeqnarray}
 where $c_i^n$ are arbitrary but fixed natural number.
From Lemma \ref{sample1}, it follows that  
\begin{IEEEeqnarray}{lCl} 
 (g_j^{(w,n)})^r&= &(g_0^{(w,n)})^*)^{s}h_j^{(w,n)} \nonumber \\ 
 &\sim& 1^{\otimes (n+ \sum_{i=1}^{n-1} \ell(w_i))}\otimes 1^{\otimes (\sigma(j)-1)}\otimes 
 (\sqrt{1-q^{2d_iN}}S^*)^r\otimes 1^{\otimes (\ell(w_n)-\sigma(j))}. \nonumber 
\end{IEEEeqnarray}
on the vector  subspace generated by standard orthonormal basis elements having $e_0$ at $\big(n+\sum_{i=1}^{n-1} \ell(w_i)+\sigma(k)\big)^{th}$ place for each $k <j$.
Employing this and  Lemma \ref{independent}, we get the claim.
\qed \\
If one counts the number of linearly independent endomorphisms given in part $(2)$ of the lemma, one can show that $GKdim(\mathcal{O}_q(G)) \geq 2\ell(n)+1$. To get the upper bound, we 
need to extend this result.
%As in the previous section, we will extend this result to incorporate all parts of the  Weyl word. Note that the polynomials $g_{k*}^{(w,n)}$'s involve variables of the type $\chi_w^n((u_j^i)^*)$. Hence we need to strengthen 
%Proposition \ref{operator}. 
\blmma \label{operator1}
%Let $\mathcal{O}_q(G)$ be one of the Hopf $*$-algebras $\mathcal{O}_q(SU(n+1))$, $\mathcal{O}_q(SP(2n))$ or  $\mathcal{O}_q(\mbox{Spin}(2n))$.
 % For $w \in W$, let  $V_w=V_{w^1}\otimes V_{w^2}\otimes \cdots \otimes V_{w^n}$ be the associated   $\mathcal{O}_q(G)$-module.
 Let $w \in W$ and   $\chi_w^n$ be the associated action of $\mathcal{O}_q(G)$ in the vector space $c_{00}(\bbz)^{\otimes n} \otimes c_{00}(\bbn)^{\otimes \ell(w)}$. 
 Then for each $1 \leq i < n$, there exist endomorphisms $P_1^i,P_{2}^i, \cdots P_{N_i}^i$ and $R_1^i,R_{2}^i, \cdots R_{N_i}^i$ in $C(G_q)$ such that 
   for all $1\leq j \leq N_i$, one has
  \begin{IEEEeqnarray}{rCl}
    P_j^i(v\otimes e_0 \otimes \cdots \otimes e_0)&=&C (\chi_{w_1w_2\cdots w^i}^{i}(u_{j}^{N_i})v)\otimes e_0 \otimes \cdots \otimes e_0 \nonumber \\
   R_j^i(v\otimes e_0 \otimes \cdots \otimes e_0)&=&C \big((\chi_{w_1w_2\cdots w^i}^{i}(u_{j}^{N_i}))^*v\big)\otimes e_0 \otimes \cdots \otimes e_0 \nonumber 
  \end{IEEEeqnarray}
 where $v\in  c_{00}(\bbz)^{\otimes n} \otimes c_{00}(\bbn)^{\otimes \sum_{j=1}^{i}\ell(w_j)}$  and $C$ is a nonzero constant. 
\elmma
\prf
Define the endomorphisms 
 \begin{IEEEeqnarray}{rCl}
 P_j^i&:=&\chi_{w}^n(u_{r_{w^{i+1}w^{i+2}\cdots w^{n}}(j)}^{N_n^i}) \nonumber \\
 R_j^i&:=&(\chi_{w}^n(u_{r_{w^{i+1}w^{i+2}\cdots w^{n}}(j)}^{N_n^i}))^* \nonumber 
\end{IEEEeqnarray}
 for $1 \leq j \leq N_i$. By applying Lemma \ref{unique path},  the claim follows immediately. \qed 

\blmma \label{sample4}
%Let $\mathcal{O}_q(G)$ be one of the Hopf $*$-algebras $\mathcal{O}_q(SU(n+1))$, $\mathcal{O}_q(SP(2n))$ or  $\mathcal{O}_q(\mbox{Spin}(2n))$.
Let $w \in W$ and   $\chi_w^n$ be the associated action of $\mathcal{O}_q(G)$ in the vector space $c_{00}(\bbz)^{\otimes n} \otimes c_{00}(\bbn)^{\otimes \ell(w)}$. 
Then for each $1 \leq i \leq n$, there exist   permutations $\sigma_i$ of $\{1,2,\cdots ,\ell(w_i)\}$ and
polynomials $g_{0}^{(w,i)}, g_{1}^{(w,i)}, g_{2}^{(w,i)} \cdots ,g_{\ell(w_i)}^{(w,i)}, 
g_{1*}^{(w,i)},  \cdots ,g_{\ell(w_i)*}^{(w,i)}$%, \cdots ,g_{0}^{(w,n)}, g_{1}^{(w,n)},\cdots ,g_{\ell(w_n)}^{(w,n)}$,  \\ $g_{1*}^{(w,n)},\cdots ,g_{\ell(w_n)*}^{(w,n)}$ 
with noncommutative 
variables $\chi_w^n(u_s^{r})$'s and $\chi_w^n((u_s^{r})^*)$'s
such that 
\begin{enumerate}
 \item
 \begin{IEEEeqnarray}{lCl} 
 \overleftarrow{\prod_{i=1}^n}(g_{1*}^{(w,i)})^{p_{\sigma_i(1)}^i}(g_{1}^{(w,n)})^{r_{\sigma_i(1)}^i}
 (g_{2*}^{(w,i)})^{p_{\sigma_i(2)}^i}(g_{2}^{(w,i)})^{r_{\sigma_i(2)}^i}\cdots (g_{\ell(w_i)*}^{(w,i)})^{p_{\sigma_i(\ell(w_i))}^i}  \nonumber \\    (g_{\ell_i}^{(w,i)})^{r_{\sigma_i(\ell_i)}^i}
   \overleftarrow{\prod_{i=1}^n} (g_{0}^{(w,i)})^{r_0^i}(e_0\otimes  \cdots \otimes e_0) 
 = C e_{r_0^n}\otimes e_{r_0^{n-1}} \otimes \cdots \otimes e_{r_0^1}  \otimes  e_{r_{1}^1-p_1^1}\otimes e_{r_{2}^1-p_2^1} \otimes \cdots  \nonumber \\
  \cdots  e_{r_{\ell(w_1)}^1-p_{\ell(w_1)}^1}\otimes \cdots \otimes e_{r_{1}^n-p_1^n}\otimes   \otimes e_{r_{2}^n-p_2^n} \otimes \cdots \otimes e_{r_{\ell(w_n)}^n-p_{\ell(w_n)}^n} \nonumber 
\end{IEEEeqnarray}
where $r_0^j \in \bbn$, $r_j^i, p_j^i \in \bbn$ and $r_j^i\geq p_j^i$ for $1 \leq j \leq \ell(w_i)$ and $1 \leq i \leq n$. 
 \item  The elements of the set
 \begin{IEEEeqnarray}{lCl}
   \Big\{ \overleftarrow{\prod_{i=1}^n}(g_{1*}^{(w,i)})^{p_{\sigma_i(1)}^i}(g_{1}^{(w,n)})^{r_{\sigma_i(1)}^i}
 (g_{2*}^{(w,i)})^{p_{\sigma_i(2)}^i}(g_{2}^{(w,i)})^{r_{\sigma_i(2)}^i}\cdots (g_{\ell_i*}^{(w,i)})^{p_{\sigma_i(\ell_i)}^i}  \nonumber \\    (g_{\ell_i}^{(w,i)})^{r_{\sigma_i(\ell_i)}^i}\overleftarrow{\prod_{i=1}^n} (g_{0}^{(w,i)})^{r_0^i} : r_0^i \in \bbn, 
  r_j^i\geq p_j^i \mbox{ for } 1 \leq j \leq \ell(w_i) \mbox{ and } 1 \leq i \leq n \Big\} \nonumber % r_j^1\geq p_j^i \in \bbn,  1 \leq j \leq \ell(w_n) \mbox{ and }1 \leq i \leq n\} \nonumber
 \end{IEEEeqnarray}
are linearly independent  endomorphisms.  
\end{enumerate}
\elmma
\prf For $0 \leq j \leq \ell(w_n)$, let  $g_j^{(w,n)}$ and  $g_{j*}^{(w,n)}$ be  the polynomials as given in Lemma \ref{sample3} and  $\sigma_n$ be the associated permutation. 
To define permutations $\sigma_i$ and polynomials  $g_j^{(w,i)}$ and $g_{j*}^{(w,i)}$ for $0 \leq j \leq \ell(w_i)$ and $1\leq i<n$   we view $w_1w_2\cdots w_i$ as an element of Weyl group of $G$ of rank $i$. Therefore  we 
can define polynomial $g_j^{(w_1w_2\cdots w_i,i)}$ and the permutation $\sigma_i$ from Proposition \ref{sample3}. 
Replace the variables $\chi_{w_1w_2\cdots w_i}^i(u_k^{N_i})$ with $P_k^i$ 
and $\chi_{w_1w_2\cdots w_i}^i((u_k^{N_i})^*)$ with $R_k^i$ for $1 \leq k \leq N_i$  in the polynomials $g_j^{(w_1w_2\cdots w_i,i)}$ and $g_{j*}^{(w_1w_2\cdots w_i,i)}$ to define the polynomial $g_{j}^{(w,i)}$  and 
$g_{j*}^{(w,i)}$ respectively for  all $0\leq j \leq \ell(w_i)$. 
Now both parts of the claim follows from Lemma  \ref{sample3} and Lemma \ref{operator1}. 
\qed  \\
Define
\begin{IEEEeqnarray}{rCl}
 \xi_{G_q}=\begin{cases}
            \{u_j^i:1\leq i,j \leq N_n\} \cup \{1\} \quad \quad  \mbox{ for } \mathcal{O}_q(G)=\mathcal{O}(SU_q{(n+1)}) \mbox{ or } \mathcal{O}_q(SP(2n)), \cr
            \{u_j^i:1\leq i,j \leq N_n\} \cup \{z_j^i:1\leq i,j \leq 2^n\} \cup \{1\}  \quad  \mbox{ for } \mathcal{O}_q(G)=\mathcal{O}(\mbox{Spin}_q(n)). \cr
           \end{cases} \nonumber
\end{IEEEeqnarray} 
Then $\xi_{G_q}$ is a generating set of $\mathcal{O}_q(G)$ containing $1$. We will now prove our main result.
%Now we have all the ingredients required to prove the main result of this section.
\bthm \label{GKalgebra}
%Let $\mathcal{O}_q(G)$ be one of the Hopf $*$-algebras $\mathcal{O}_q(SU(n+1))$, $\mathcal{O}_q(SP(2n))$ or  $\mathcal{O}_q(\mbox{Spin}(2n))$. 
Let $\omega_n$ be  the longest element of the Weyl group  of $G$. Then  one has 
\[
 \mbox{GKdim}(\mathcal{O}_q(G))= 2\ell(\omega_n)+n.
\]
\ethm
\prf From Lemma \ref{lessthandim}, it is enough  to show that $\mbox{GKdim}(\mathcal{O}_q(G))\geq 2\ell(w_n)+n$. Since the homomorphism $\chi_{\omega_n}^n$ is faithful, we will
without loss of generality work with the algebra $\chi_{\omega_n}^n(\mathcal{O}_q(G))$.
Take the generating set $\varUpsilon$ to be $\chi_{\omega_n}^n(\xi_{G_q}) \cup \chi_{\omega_n}^n(\xi_{G_q}^*)$. Define 
\[
 M_0:=\max\{\mbox{degree of } g_j^{(\omega_n,i)}:  0\leq j\leq \ell((\omega_n)_i), 1\leq i\leq n\}.
\]
 Then by part (2) of Lemma \ref{sample4}, we have
\begin{IEEEeqnarray}{rCl}
 \mbox{GKdim}(\mathcal{O}_q(G)) &\geq& \varlimsup \frac{\ln \dim(\varUpsilon^{M_0k})}{\ln M_0k} 
 \geq \varlimsup\frac{\ln\frac{{k+2\ell(\omega_n)+n-1\choose k}}{2}}{\ln M_0k}  
 =2\ell(\omega_n)+n. \nonumber 
\end{IEEEeqnarray}
This completes the proof.	
\qed 

 \brmrk
 The proof we have given here is very rigid in the sense that it largely depends 
 upon a particular way of representing the algebra. There must be   a canonical way of computing GKdim of these algebras.
In our view, the main obstruction  is to get Lemma \ref{unique path} in  a more general set up.
 \ermrk

\bcrlre  \label{GKalgebracor}
One has
\begin{itemize}
 \item 
 GKdim$(\mathcal{O}_q(SU(n+1)))=n^2+2n=$dim$(SU(n+1))$.
 \item 
 GKdim$(\mathcal{O}_q(SP(2n)))=2n^2+n=$dim$(SP(2n))$.
 \item 
GKdim$(\mathcal{O}_q(SO(2n)))=\mbox{GKdim}(\mathcal{O}_q(\mbox{Spin}(2n)))=2n^2-n=$dim$(SO(2n))$.
%\item 
% GKdim$(C[SU_q(n+1)/SU_q(n+1-m)])=m(2n-m+2)=$dim$(SU(n+1)/SU(n+1-m))$.
% \item 
 %GKdim$(C[SP_q(2n)/SP_q(2n-2m)])=m(4n-2m+1)=$dim$(SP(2n)/SP(2n-2m))$.
\end{itemize}
\ecrlre
\prf Let $\omega_n$ be  the longest word of $\mathcal{O}_q(G)$. \\
\textbf{Case 1:} $\mathcal{O}_q(SU(n+1))$.\\ 
In this case, the $r^{th}$-part  of $\omega_n$ is $s_rs_{r-1} \cdots s_1$ for $1 \leq r \leq n$.  Hence $\ell(\omega_n)=n(n+1)/2$. Therefore by 
Theorem \ref{GKalgebra}, we have 
\[
 \mbox{GKdim}(\mathcal{O}_q(SU(n+1)))=\frac{2n(n+1)}{2}+n=n^2+2n.
\]
\textbf{Case 2:} $\mathcal{O}_q(SP(2n))$.\\ 
For each  $1 \leq r \leq n$, the $r^{th}$-part of $\omega_n$ is $s_rs_{r+1}\cdots s_{n-1}s_ns_{n-1}\cdots s_{r+1}s_r$.  Hence by applying Theorem \ref{GKalgebra}, we get
\[
 \mbox{GKdim}(\mathcal{O}_q(SP(2n)))=2\ell(\omega_n)+n=2n^2+n.
\]
\textbf{Case 3:} $\mathcal{O}_q(\mbox{Spin}(2n))$.\\ 
For each  $1 \leq r \leq n$, the $r^{th}$-part of $\omega_n$ is $s_rs_{r+1}\cdots s_{n-2}s_{n-1}s_ns_{n-2}\cdots s_{r+1}s_r$. Hence $\ell(\omega_n)=n^2-n$. Therefore from Theorem \ref{GKalgebra}, we get
\[
 \mbox{GKdim}(\mathcal{O}_q(\mbox{Spin}(2n)))=2\ell(\omega_n)+n=2n^2-n.
\]
Moreover  since the polynomials $g_0^{(\omega_n,k)}$,  $g_j^{(\omega_n,k)}$ and $g_{j*}^{(\omega_n,k)}$ %for $1 \leq j \leq \ell_k$ and $n-m+1 \leq k \leq n$  
given 
in Lemma \ref{sample4} involve variables $u_j^i$'s which are in $\mathcal{O}(SO_q(2n)$, we get 
\[
\mbox{GKdim}(\mathcal{O}_q(SO(2n)))=\mbox{GKdim}(\mathcal{O}_q(\mbox{Spin}(2n))). 
\]
 This completes the proof. 
\qed

\newsection{Quotient spaces}
Fix a subset $S \subset \Pi$ and a subgroup $L$ of $\bbbt^{\# S^{c}}$. Let $\mathcal{O}(G_q/K_{q}^{S,L})$ be the quotient Hopf $*$-subalgebra  of $\mathcal{O}_q(G)$ (see page 5, \cite{NesTus-2012ab}).
If $S$ is the empty set $\phi$, define $W_{\phi}=\{id\}$. For a nonempty set $S$, define $W_S$ to be the subgroup of $W_n$ generated by the  simple reflections $s_{\alpha}$ with $\alpha \in S$. 
%Let $W^S$ be the set of elements $w$ such that $w(\alpha) >0$ for all $\alpha \in S$.
Let
\[
 W^S:=\{w \in W_n : \ell(s_{\alpha}w)>\ell(w) \quad \forall \alpha \in S\}.
\]
Define the algebra $\mathcal{P}(C(L))$ to be the quotient of $\mathcal{P}(C(\bbbt^m))$ by the ideal consisting of polynomials vanishing in  $L$.
\bthm \label{faithful1}
Let $\omega_{n}^S$ be the longest word of $W^S$, $m$ be the cardinality of $S$  and $k$ be the rank of $L$. 
Then the homomorphism
\[
                            \chi_{\omega_{n}^{S}}^n : \mathcal{O}(G_q/K_{q}^{S,L})\longrightarrow \mathcal{P}(C(\bbbt))^{\otimes m}\otimes \mathcal{P}_q(\scrt)^{\otimes \ell(\omega_n^{S}})
\]
is faithful. Moreover, the image  $\chi_{\omega_{n}^{S}}^n ( \mathcal{O}(G_q/K_{q}^{S,L})) $ is 
contained in the algebra $\mathcal{P}(C(L))\otimes \mathcal{P}_q(\scrt)^{\otimes \ell(\omega_n^{S}})$.
\ethm

\blmma \label{iso}
Let $L$ be a subgroup of $\bbbt^m$ of rank $k$. Then 
there exists an algebra isomorphism 
\[
 \Phi:\mathcal{P}(C(L)) \rightarrow \mathcal{P}(C(\bbbt))^{\otimes k}.
\]
\elmma
\prf  Using Fourier transform, one can identify dual group $\hat{L}$ as a subgroup of $\bbz^m$ isomorphic to $\bbz^k$. Fix a linear isomorphism $\phi: \bbz^k \rightarrow \hat{L}$. 
Applying inverse Fourier transform and using Pontriagin duality, one can identify  points of $\bbbt^k$ with points of $L$  via the monomial map $\hat{\phi}$. This induces an isomorphism 
\[
 \Phi:\mathcal{P}(C(L)) \rightarrow \mathcal{P}(C(\bbbt))^{\otimes k}
\]
such that $\Phi(g)(t_1,t_2,\cdots ,t_k)=g(\hat{\phi}(t_1,t_2,\cdots,t_k))$.
\qed
\brmrk \label{rem}
Suppose that the polynomials $\{g_i(t_1,t_2,\cdots t_m):1 \leq i \leq m\}$ generate the algebra $\mathcal{P}(C(L))$. Then  from 
the Lemma \ref{iso}, there exist monomials  $\{h_i(g_1,g_2,\cdots g_m):1 \leq i \leq k\}$ such that $h_i(g_1,g_2,\cdots g_m)(t_1,t_2,\cdots t_m)=t_i$ for $1 \leq i \leq k$. Hence the elements of the set
\[
 \{\prod_{i=1}^n(h_i(g_1,g_2,\cdots g_m))^{r_i}: r_i \in \bbn, 1 \leq i \leq k \}
\]
are linearly independent  polynomials.
\ermrk
Let $S_1$ be the empty subset of $\Pi$.
For $2 \leq m \leq n$,   define $S_m$ to be the set  $\{1,2,\cdots ,m-1\}$ if $\mathcal{O}_q(G)=\mathcal{O}(SU_q{(n+1)})$ and $\{n-m+2, \cdots ,n\}$ if $\mathcal{O}_q(G)= \mathcal{O}_q(SP(2n))$.
If $L=\bbbt^{m}$ then $\mathcal{O}(G_q/K_{q}^{S_{n-m+1},L})$ is same as $\mathcal{O}(SU_q(n+1)/SU_q(n+1-m))$ or $\mathcal{O}(SP_q(2n)/SP_q(2n-2m))$  
if $\mathcal{O}_q(G)$ is $\mathcal{O}_q(SU(n+1))$ or $\mathcal{O}_q(SP(2n))$  respectively.
\bthm \label{GKalgebra1} 
Let $\omega_n^{S_{n-m+1}}$ be  the longest element of the $W^{S_{n-m+1}}$ and $k$ be the rank of $L$. Then  one has 
\[
 \mbox{GKdim }\mathcal{O}(G_q/K_{q}^{S_{n-m+1},L})= 2\ell(\omega_n^{S_{n-m+1}})+k.
\]
\ethm
\prf By Theorem \ref{faithful1}, the algebra  $\mathcal{O}(G_q/K_{q}^{S_{n-m+1},L})$ can be viewed as a subalgebra of $\mathcal{P}(C(\bbbt))^{\otimes k}\otimes \mathcal{P}_q(\scrt)^{\otimes \ell(\omega_n^{S_{n-m+1}}})$. 
Hence using properties of GKdim, we get 
\[
 \mbox{GKdim }\mathcal{O}(G_q/K_{q}^{S_{n-m+1},L})\leq 2\ell(\omega_n^{S_{n-m+1}})+k.
\]

To show the equality,  observe that 
\begin{enumerate}
 \item 
 for an element $w \in W^{S_{n-m+1}}$  and $1 \leq r \leq n-m$, the $r^{th}$-part $w_r=\psi_{r,k_{r}}^{(\epsilon_{r})}(w)$ is identity element of $W_n$. Hence $w$ can be written uniquely as $w= w_{n-m+1} w_{n-m+2}\cdots w_{n}$.
 \item 
 It follows from the definition  that  entries of last 
$m$ rows of $(\!(u_j^i)\!)$ is in the quotient algebra $\mathcal{O}(G_q/K_{q}^{S{n-m+1},L})$.
 \item
 The polynomials $g_j^{(w,k)}$  and $g_{j*}^{(w,k)}$for $1 \leq j \leq \ell_k$ and $n-m+1 \leq k \leq n$ involve variables consisting of entries of last 
$m$ rows of $(\!(u_j^i)\!)$.
\item 
It follows from equation (\ref{eq1}) that for $1 \leq i \leq m$, we have 
\[
 g_{0}^{(w,i)} \sim g_i(t_1,t_2,\cdots, t_m) \otimes 1^{\otimes \ell(w)}
\]
where  $g_i \in \mathcal{P}(C(L))$ is the projection function $t_i$ on $i^{th}$ co-ordinate  restricted to $L$. Moreover, $\{g_i(t_1,t_2,\cdots t_m):1 \leq i \leq m\}$ generate the algebra $\mathcal{P}(C(L))$. 
\item  Using remark (\ref{rem}) and Lemma \ref{sample4}, one can show that the elements of the set
 \begin{IEEEeqnarray}{lCl}
   \Big\{ \overleftarrow{\prod_{i=n-m+1}^n}(g_{1*}^{(w,i)})^{p_{\sigma_i(1)}^i}(g_{1}^{(w,n)})^{r_{\sigma_i(1)}^i}
 (g_{2*}^{(w,i)})^{p_{\sigma_i(2)}^i}(g_{2}^{(w,i)})^{r_{\sigma_i(2)}^i}\cdots (g_{\ell(w_i)*}^{(w,i)})^{p_{\sigma_i(\ell(w_i))}^i}  \nonumber \\    
 (g_{\ell(w_i)}^{(w,i)})^{r_{\sigma_i(\ell(w_i))}^i}\overleftarrow{\prod_{i=1}^n} \big(h_i(g_{0}^{(w,1)},g_{0}^{(w,2)}, \cdots ,g_{0}^{(w,k)})^{r_0^i} : r_0^i \in \bbn, 
  r_j^i\geq p_j^i \nonumber \\  \mbox{ for } 1 \leq j \leq \ell(w_i) \mbox{ and } n-m+1 \leq i \leq n \Big\} \nonumber % r_j^1\geq p_j^i \in \bbn,  1 \leq j \leq \ell(w_n) \mbox{ and }1 \leq i \leq n\} \nonumber
 \end{IEEEeqnarray}
are linearly independent  endomorphisms. 	      
\end{enumerate}
With these facts,  the same arguments  used in Theorem \ref{GKalgebra} will prove the claim.                                                                                                              
\qed

\bcrlre
One has
\begin{itemize}
 \item 
 GKdim$(\mathcal{O}(SU_q(n+1)/SU_q(n+1-m)))=$dim$(SU(n+1)/SU(n+1-m))$.
 \item 
 GKdim$(\mathcal{O}(SP_q(2n)/SP_q(2n-2m)))=$dim$(SP(2n)/SP(2n-2m))$.
\end{itemize}
\ecrlre
\prf Let $\omega_n^{S_{n-m+1}}$ be  the longest word of $\mathcal{O}_q(G)$. \\
\textbf{Case 1:} $\mathcal{O}(SU_q(n+1)/SU_q(n+1-m))$.\\ 
In this case, the $r^{th}$-part  of $\omega_n^{S_{n-m+1}}$ is  % identity for $1 \leq r \leq n-m$ and 
$s_rs_{r-1} \cdots s_1$ for $n-m+1 \leq r \leq n$.  Hence $\ell(\omega_n^{S_{n-m+1}})=\frac{n(n+1)-(n-m)(n-m+1)}{2}$. Therefore by 
Theorem \ref{GKalgebra}, we have 
\begin{IEEEeqnarray}{lCl}
  \mbox{GKdim}(\mathcal{O}_q(SU(n+1)))&=&n(n+1)-(n-m)(n-m+1)+m \nonumber\\
  &=& n(n+1)+n-(n-m)(n-m+1)+(n-m) \nonumber \\
  &=& \mbox{dim}(SU(n+1)) -\mbox{dim}(SU(n+1-m))  \mbox{ (by Corollary } \ref{GKalgebracor}) \nonumber\\
  &=&\mbox{dim}(SU(n+1)/SU(n+1-m)). \nonumber
\end{IEEEeqnarray}
\textbf{Case 2:} $\mathcal{O}(SP_q(2n)/SP_q(2n-2m))$.\\ 
In this case, the $r^{th}$-part  of $\omega_n^{S_{n-m+1}}$ is % identity for $1 \leq r \leq n-m$ and 
$s_rs_{r+1}\cdots s_{n-1}s_ns_{n-1}\cdots s_{r+1}s_r$ for $n-m+1 \leq r \leq n$.  Hence by applying Theorem \ref{GKalgebra} and following the steps of part $(1)$, we get the claim.
\qed

\noindent{\sc Partha Sarathi Chakraborty} (\texttt{parthac@imsc.res.in})\\
         {\footnotesize Institute of Mathematical Sciences (HBNI),   CIT Campus,
Taramani,
Chennai, 600113, INDIA}

\noindent{\sc Bipul Saurabh} (\texttt{saurabhbipul2@gmail.com})\\
         {\footnotesize Institute of Mathematical Sciences (HBNI), CIT Campus,
Taramani,
Chennai, 600113, INDIA}

\end{document}